\documentclass[oneside, reqno] {amsart}

\usepackage{xypic}
\input xy
\xyoption{all}
\usepackage{epsfig}
\usepackage{amsthm}
\usepackage{amssymb}
\usepackage{amsmath}
\usepackage{amscd}
\usepackage{color}
\usepackage[T1]{fontenc}
\usepackage{upgreek}
\usepackage[font=scriptsize]{caption}
\usepackage{wrapfig}

\usepackage{graphicx}
\usepackage{subfigure}

\newcommand{\bg}{\begin{equation}}
\newcommand{\ed}{\end{equation}}
\newcommand{\bga}{\begin{eqnarray}}
\newcommand{\eda}{\end{eqnarray}}

\def\cbdu{\par{\raggedleft$\Box$\par}}

\newtheorem {Theorem}  {Theorem}

\numberwithin{Theorem}{section}

\newtheorem {Lemma}[Theorem]  {Lemma}

\theoremstyle{definition}

\theoremstyle{remark}

\expandafter\chardef\csname pre amssym.def
at\endcsname=\the\catcode`\@ \catcode`\@=11
\def\undefine#1{\let#1\undefined}
\def\newsymbol#1#2#3#4#5{\let\next@\relax
 \ifnum#2=\@ne\let\next@\msafam@\else
 \ifnum#2=\tw@\let\next@\msbfam@\fi\fi
 \mathchardef#1="#3\next@#4#5}
\def\mathhexbox@#1#2#3{\relax
 \ifmmode\mathpalette{}{\m@th\mathchar"#1#2#3}%
 \else\leavevmode\hbox{$\m@th\mathchar"#1#2#3$}\fi}
\def\hexnumber@#1{\ifcase#1 0\or 1\or 2\or 3\or 4\or 5\or 6\or 7\or 8\or
 9\or A\or B\or C\or D\or E\or F\fi}

\font\teneufm=eufm10 \font\seveneufm=eufm7 \font\fiveeufm=eufm5
\newfam\eufmfam
\textfont\eufmfam=\teneufm \scriptfont\eufmfam=\seveneufm
\scriptscriptfont\eufmfam=\fiveeufm

\catcode`\@=\csname pre amssym.def at\endcsname

\newcounter{remark}
\def  \12  {{\frac{1}{2}}}

\def\build#1_#2^#3{\mathrel{\mathop{\kern 0pt#1}\limits_{#2}^{#3}}}

\numberwithin{equation}{section}

\begin{document}

\title[EMHD well-posedness]{Well-posedness of the electron MHD with partial resistivity}


\author [Mimi Dai]{Mimi Dai}

\address{Department of Mathematics, Statistics and Computer Science, University of Illinois at Chicago, Chicago, IL 60607, USA}
\email{mdai@uic.edu} 

\author [Hassan Babaei]{Hassan Babaei}

\address{Department of Mathematics, Statistics and Computer Science, University of Illinois at Chicago, Chicago, IL 60607, USA}
\email{hbabae2@uic.edu} 

\thanks{The authors are partially supported by the NSF grant DMS--2308208 and Simons Foundation.}

\begin{abstract}
Due to the singular nonlinear Hall term, the non-resistive electron magnetohydrodynamics (MHD) is not known to be locally well-posed in general. In this paper we consider the $2\frac12$D electron MHD with either horizontal or vertical resistivity and show local well-posedness in Sobolev spaces.

\bigskip


\hspace{0.02cm}CLASSIFICATION CODE: 35Q35, 76B03, 76D09, 76E25, 76W05.
\end{abstract}

\maketitle

\section{Introduction}



The electron magnetohydrodynamics (MHD) system 
\begin{equation}\label{emhd}
\begin{split}
B_t+ \nabla\times ((\nabla\times B)\times B)=&\ \nu\Delta B,\\
\nabla\cdot B=&\ 0
\end{split}
\end{equation}
with $B$ being the unknown magnetic field and $\nu\geq 0$ the resistivity parameter, is an important model to capture rapid magnetic reconnection phenomena in plasma physics. As demonstrated in our previous work \cite{Dai-emhd-25}, the Hall term $\nabla\times ((\nabla\times B)\times B)$ has an intricate nonlinear structure which tends to show singular behavior. The local well-posedness of the non-resistive electron MHD \eqref{emhd}, i.e. $\nu=0$, remains a major open problem from the perspective of mathematics. For relevant topics on well-posedness, ill-posedness and regularity of the electron MHD, we refer the reader to the incomplete list of works \cite{Dai-emhd-24, Dai-emhd-18, Dai-Oh, JO1, JO2, JO3} and references within these papers.

A special context of the electron MHD in $2\frac12$D setting was initially (mathematically) studied in \cite{Dai-W}. Consider the magnetic field $B$ independent of the $z$ variable. It can be written as $B=\nabla\times (a\vec e_z)+b\vec e_z$ with $a=a(x,y,t)$ and $b=b(x,y,t)$. It follows from \eqref{emhd} that
\begin{equation}\label{eq-ab}
\begin{split}
\partial_t a+a_yb_x-a_xb_y&=\mu\Delta a,\\
\partial_t b-a_y\Delta a_x+a_x\Delta a_y&=\nu \Delta b
\end{split}
\end{equation}
with the assumption that the horizontal and vertical resistivity parameters $\mu$ and $\nu$ are not necessarily the same. 


In our previous work \cite{Dai-emhd-23}, system \eqref{eq-ab} with only vertical resistivity, i.e. $\mu=0$ and $\nu>0$, was shown to have a global solution in high regularity Sobolev space near a steady state $(a^*, b^*)=(y,0)$ of the shear type. On the other hand, with $\mu=\nu=0$, we \cite{Dai-emhd-24} proved that system \eqref{eq-ab} is ill-posed in the sense of norm inflation in Sobolev space $(a,b)\in\dot H^{\beta}\times \dot H^{\beta-1}$ with $3<\beta<4$.

In this paper we aim to show local well-posedness of \eqref{eq-ab} with either vertical or horizontal resistivity. 

\begin{Theorem}\label{thm}
Let $(a_0, b_0)\in H^{s+1}(\mathbb R^2)\times H^{s}(\mathbb R^2)$ be the initial data with $s>2$. Assume either $\{\mu>0, \nu=0\}$ or $\{\mu=0, \nu>0\}$. There exists $T>0$ depending on $\|a_0\|_{H^{s+1}}$, $\|b_0\|_{H^{s}}$ and $\mu$ (or $\nu$) such that system \eqref{eq-ab} has a unique solution $(a(t), b(t))\in H^{s+1}(\mathbb R^2)\times H^{s}(\mathbb R^2)$ on $[0,T)$ with $(a(t), b(t))=(a_0,b_0)$. 
\end{Theorem}

\bigskip

\section{Preliminaries}
\label{sec-pre}

\subsection{Notations}
Through the rest of the text, $C$ is used as an absolute general constant. The inequality $A\lesssim B$ means $A\leq CB$ for some constant $C>0$. Denote the variable $\vec x=(x,y)$.

\subsection{Basic energy estimate}
\label{sec-basic}
System \eqref{eq-ab} satisfies a basic energy estimate. Formally, multiplying the first equation of \eqref{eq-ab} by $\Delta a$ and the second one by $b$, it follows from integrating and adding the two equations
\begin{equation}\notag
\begin{split}
&\frac12\frac{d}{dt}\int_{\mathbb R^2}|\nabla a|^2+|b|^2\,dxdy+\int_{\mathbb R^2}\mu|\Delta a|^2+\nu|\nabla b|^2\,dxdy\\
=& \int_{\mathbb R^2}\nabla^\perp b\cdot\nabla a\Delta a\, dxdy+ \int_{\mathbb R^2}\nabla^\perp \Delta a\cdot\nabla a b\, dxdy\\
=&-  \int_{\mathbb R^2} b\nabla a\cdot \nabla^\perp\Delta a\, dxdy+ \int_{\mathbb R^2}\nabla^\perp \Delta a\cdot\nabla a b\, dxdy\\
=&\ 0
\end{split}
\end{equation}
where we used integration by parts and the fact $\nabla^\perp\cdot \nabla a=0$. Hence we have the a priori estimates 
\[\nabla a, b\in L^\infty_t L^2_x.\]

Moreover, due to the transport feature of the $a$ equation of \eqref{eq-ab}, it is clear that if $a_0\in L^p_x$ for any $1\leq p\leq \infty$, then $a\in L^\infty_t L^p_x$.

\subsection{Littlewood-Paley projection and commutator estimate}

Denote the frequency number $\lambda_q=2^q$ for any integer $q\geq -1$. Let $\Delta_q u$ be the $q$-th Littlewood-Paley projection of $u$. We adapt the simplified notations
\[\Delta_q u=u_q, \ \ \ \ u_{\leq Q}=\sum_{q\geq-1}^Q u_q, \ \ \ \ \widetilde u=u_{q-1}+u_q+u_{q+1}.\]
It is worth noting the Sobolev norm $\|u\|_{\dot H^s}$ is equivalent to 
\[\left(\sum_{q=-1}^\infty \lambda_q^{2s} \|u_q\|_{L^2}^2\right)^{\frac12}.\]

We define the commutator
\[[\Delta_q, f] g=\Delta_q(fg)-f g_q.\]

\begin{Lemma}\label{le-comm} \cite{Dai-emhd-1d}
For any $1<p<\infty$, the estimate
\begin{equation}\notag
\left\|[\Delta_q, f_{\leq p-2}] \nabla g_q\right\|_{L^p}\leq C \|\nabla f_{\leq p-2}\|_{L^\infty}\|g_q\|_{L^p}
\end{equation} 
holds for a constant $C$ depending on $p$.
\end{Lemma}

\bigskip

\section{Proof of Theorem \ref{thm} for $\mu>0$ and $\nu=0$}
\label{sec-proof1}

In this section we prove the local well-posedness of Theorem \ref{thm} in the case of only horizontal resistivity, i.e. $\mu>0$ and $\nu=0$ in \eqref{eq-ab}. As customary, the key ingredient is to establish a priori estimate $(a(t),b(t))\in H^{s+1}\times H^s$ for $s>2$. The existence and uniqueness of the solution follow from the a priori estimate by standard argument. 

Taking the $q$-th Littlewood-Paley projection $\Delta_q$ on the first equation of \eqref{eq-ab}
\[a_t=\mu\Delta a-(a_yb_x-a_xb_y)\]
and multiplying by $\lambda_q^{2s}\Delta a_{q}$, it then follows by integrating and summation over $q\geq-1$
\begin{equation}\label{priori-1}
\begin{split}
&\frac{1}{2}\frac{d}{dt} \sum_{q\geq-1} \lambda_q^{2s}\int_{\mathbb R^2} |\nabla a_q|^2\, d\vec x+\mu \sum_{q\geq-1} \lambda_q^{2s}\int_{\mathbb R^2} (\Delta a_q)^2\, d\vec x\\
=& \sum_{q\geq-1} \lambda_q^{2s}\int_{\mathbb R^2} (a_yb_x)_q(a_{xx})_{q}\, d\vec x+\sum_{q\geq-1} \lambda_q^{2s}\int_{\mathbb R^2} (a_yb_x)_q(a_{yy})_{q}\, d\vec x\\
&- \sum_{q\geq-1} \lambda_q^{2s}\int_{\mathbb R^2} (a_xb_y)_q(a_{xx})_{q}\, d\vec x-\sum_{q\geq-1} \lambda_q^{2s}\int_{\mathbb R^2} (a_xb_y)_q(a_{yy})_{q}\, d\vec x\\
=&: I_1+I_2+I_3+I_4.
\end{split}
\end{equation}
Applying similar operations to the second equation of \eqref{eq-ab} 
\[ b_t=\nu\Delta b+(a_y\Delta a_x-a_x\Delta a_y)\]
and multiplying $\lambda_q^{2s}b_{q}$ gives
\begin{equation}\label{priori-2}
\begin{split}
&\frac{1}{2}\frac{d}{dt} \sum_{q\geq-1} \lambda_q^{2s}\int_{\mathbb R^2} b_q^2\, d\vec x+ \nu\sum_{q\geq-1} \lambda_q^{2s}\int_{\mathbb R^2} |\nabla b_q|^2\, d\vec x\\
=& -\sum_{q\geq-1} \lambda_q^{2s}\int_{\mathbb R^2} a_ya_{xx}(b_{x})_{qq}\, d\vec x
-\sum_{q\geq-1} \lambda_q^{2s}\int_{\mathbb R^2} (a_ya_{yy})_q(b_{x})_{q}\, d\vec x\\
&+ \sum_{q\geq-1} \lambda_q^{2s}\int_{\mathbb R^2} (a_xa_{xx})_q(b_y)_{q}\, d\vec x 
+\sum_{q\geq-1} \lambda_q^{2s}\int_{\mathbb R^2} (a_xa_{yy})_q(b_{y})_{q}\, d\vec x\\
=&: K_1+K_2+K_3+K_4
\end{split}
\end{equation}
where we discovered cancellations to arrive at the last step. 

Adding \eqref{priori-1} and \eqref{priori-2} leads to the energy identity
\begin{equation}\label{priori-3}
\begin{split}
&\frac{1}{2}\frac{d}{dt} \sum_{q\geq-1} \lambda_q^{2s}\left(\|\nabla a_q\|_{L^2}^2+\|b_q\|_{L^2}^2\right)\\
&+\mu \sum_{q\geq-1} \lambda_q^{2s}\|\Delta a_q\|_{L^2}^2+\nu \sum_{q\geq-1} \lambda_q^{2s}\|\nabla b_q\|_{L^2}^2\\
=&: I_1+I_2+I_3+I_4+K_1+K_2+K_3+K_4.
\end{split}
\end{equation}
To obtain optimal estimates, we need to explore the cancellations in $I_1+K_1$, $I_2+K_2$, $I_3+K_3$ and $I_4+K_4$.  On the other hand, we observe it is sufficient to show details for the estimate of $I_1+K_1$ and $I_2+K_2$. Indeed, commuting the variables $x$ and $y$ in $I_3+K_3$ and $I_4+K_4$, that is, defining
\[\bar I_3(x,y)+\bar K_3(x,y)= I_3(y,x)+K_3(y,x), \]
\[\bar I_4(x,y)+\bar K_4(x,y)= I_4(y,x)+K_4(y,x).\]
We see that 
\[\bar I_3(x,y)+\bar K_3(x,y)= I_2(x,y)+K_2(x,y), \]
\[ \bar I_4(x,y)+\bar K_4(x,y)= I_1(x,y)+K_1(x,y).\]
Therefore we only provide details to estimate $I_1+K_1$ and $I_2+K_2$ in the following.

\subsection{Estimate of $I_1+K_1$}
\label{sec-I1}
We frist decompose $I_1$ and $K_1$ using Bony's paraproduct 
\begin{equation}\notag
\begin{split}
I_1
=&\sum_{q\geq -1}\sum_{|p-q|\leq 2}\lambda_q^{2s}\int_{\mathbb T^2}(a_{y,\leq p-2}b_{x,p})_q a_{xx,q} \, d\vec x\\
&+\sum_{q\geq -1}\sum_{|p-q|\leq 2}\lambda_q^{2s}\int_{\mathbb T^2}(b_{x,\leq p-2}a_{y,p})_q a_{xx,q} \, d\vec x\\
&+\sum_{q\geq -1}\sum_{p\geq q- 2}\lambda_q^{2s}\int_{\mathbb T^2}(\widetilde a_{y,p}b_{x,p})_q a_{xx,q} \, d\vec x\\
=&: I_{11}+I_{12}+I_{13},
\end{split}
\end{equation}
\begin{equation}\notag
\begin{split}
K_1
=&-\sum_{q\geq-1}\sum_{|p-q|\leq 2} \lambda_q^{2s}\int_{\mathbb T^2} (a_{y,\leq p-2}a_{xx, p})_qb_{x,q}\, d\vec x\\
&-\sum_{q\geq-1}\sum_{|p-q|\leq 2} \lambda_q^{2s}\int_{\mathbb T^2} (a_{xx,\leq p-2}a_{y, p})_qb_{x,q}\, d\vec x\\
&-\sum_{q\geq-1}\sum_{p\geq q- 2} \lambda_q^{2s}\int_{\mathbb T^2} (\widetilde a_{y,p}a_{xx, p})_qb_{x,q}\, d\vec x\\
=&: K_{11}+K_{12}+K_{13}.
\end{split}
\end{equation}
Applying the commutators 
\begin{equation}\notag
\begin{split}
[\Delta_q, a_{y,\leq p-2}\partial_x] b_p=&\ \Delta_q\left(a_{y,\leq p-2}b_{x,p}\right)-a_{y,\leq p-2}(b_{x,p})_q,\\
[\Delta_q, a_{y,\leq p-2}\partial_x] a_{x,p}=&\ \Delta_q\left(a_{y,\leq p-2}a_{xx,p}\right)-a_{y,\leq p-2}(a_{xx,p})_q,
\end{split}
\end{equation}
we further decompose $I_{11}$ and $K_{11}$ as
\begin{equation}\notag
\begin{split}
I_{11}=&\sum_{q\geq -1}\sum_{|p-q|\leq 2}\lambda_q^{2s} \int_{\mathbb T^2} [\Delta_q, a_{y,\leq p-2}\partial_x] b_p a_{xx,q} \, d\vec x\\
&+\sum_{q\geq -1}\sum_{|p-q|\leq 2}\lambda_q^{2s} \int_{\mathbb T^2} a_{y,\leq q-2}(b_{x,p})_q a_{xx,q} \, d\vec x\\
&+\sum_{q\geq -1}\sum_{|p-q|\leq 2}\lambda_q^{2s} \int_{\mathbb T^2} \left(a_{y,\leq q-2}-a_{y,\leq p-2}\right)(b_{x,p})_q a_{xx,q} \, d\vec x\\
=&:I_{111}+I_{112}+I_{113},
\end{split}
\end{equation}
\begin{equation}\notag
\begin{split}
K_{11}=&-\sum_{q\geq -1}\sum_{|p-q|\leq 2}\lambda_q^{2s} \int_{\mathbb T^2} [\Delta_q, a_{y,\leq p-2}\partial_x] a_{x,p} b_{x,q} \, d\vec x\\
&-\sum_{q\geq -1}\sum_{|p-q|\leq 2}\lambda_q^{2s} \int_{\mathbb T^2} a_{y,\leq q-2}(a_{xx,p})_q b_{x,q} \, d\vec x\\
&-\sum_{q\geq -1}\sum_{|p-q|\leq 2}\lambda_q^{2s} \int_{\mathbb T^2} \left(a_{y,\leq q-2}-a_{y,\leq p-2}\right)(a_{xx,p})_q b_{x,q} \, d\vec x\\
=&:K_{111}+K_{112}+K_{113}.
\end{split}
\end{equation}
Since for any fixed $q\geq -1$, we have
 \[\sum_{|p-q|\leq 2}(b_{x,p})_q=b_{x,q}, \ \ \ \sum_{|p-q|\leq 2}(a_{xx,p})_q=a_{xx,q},\]
 it leads to the cancellation 
\[I_{112}+K_{112}=0.\]

We estimate the rest of the terms in $I_{11}+K_{11}$ in the following. 
Applying H\"older's inequality and Young's inequality we deduce
\begin{equation}\notag
\begin{split}
|I_{111}|&\leq \sum_{q\geq-1}\sum_{|p-q|\leq 2} \lambda_q^{2s} \|a_{xy, \leq p-2}\|_{L^\infty}\|b_p\|_{L^2}\|a_{xx,q}\|_{L^2}\\
&\lesssim \sum_{q\geq-1} \lambda_q^{2s+2}\|a_{q}\|_{L^2}\|b_q\|_{L^2}\sum_{p\leq q}\lambda_p^{3} \|a_{p}\|_{L^2}\\
&\lesssim \sum_{q\geq-1} \left(\lambda_q^{s+2}\|a_q\|_{L^2}\right) \left(\lambda_q^{s}\|b_{q}\|_{L^2}\right)\sum_{p\leq q}\lambda_p^{s+1} \|a_{p}\|_{L^2}\lambda_{p}^{2-s}\\
&\leq \frac{\mu}{64}\sum_{q\geq-1} \lambda_q^{2s+4}\|a_{q}\|_{L^2}^2+\frac{C}{\mu} \sum_{q\geq-1} \lambda_q^{2s}\|b_q\|_{L^2}^2\left(\sum_{p\leq q}\lambda_p^{s+1} \|a_{p}\|_{L^2}\lambda_{p}^{2-s}\right)^2.
\end{split}
\end{equation}
Since $s>2$ we have, using H\"older's inequality for sequence gives
\begin{equation}\notag
\begin{split}
|I_{111}|
&\leq \frac{\mu}{64}\sum_{q\geq-1} \lambda_q^{2s+4}\|a_{q}\|_{L^2}^2+\frac{C}{\mu}\left(\sum_{p<\infty}\lambda_p^{s+1} \|a_{p}\|_{L^2}\lambda_{p}^{2-s}\right)^2 \sum_{q\geq-1} \lambda_q^{2s}\|b_q\|_{L^2}^2\\
&\leq \frac{\mu}{64}\sum_{q\geq-1} \lambda_q^{2s+4}\|a_{q}\|_{L^2}^2+\frac{C}{\mu}\sum_{p<\infty}\lambda_{p}^{2(2-s)}\sum_{p<\infty}\lambda_p^{2s+2} \|a_{p}\|_{L^2}^2\sum_{q\geq-1} \lambda_q^{2s}\|b_q\|_{L^2}^2\\
&\leq \frac{\mu}{64}\sum_{q\geq-1} \lambda_q^{2s+4}\|a_{q}\|_{L^2}^2+\frac{C}{\mu} \sum_{q\geq-1} \lambda_q^{2s+2}\|a_q\|_{L^2}^2\sum_{q\geq-1} \lambda_q^{2s}\|b_q\|_{L^2}^2.
\end{split}
\end{equation}
Similarly for $s>2$ we have
\begin{equation}\notag
\begin{split}
|I_{113}|&\lesssim \sum_{q\geq -1}\lambda_q^{2s}\|a_{y,q}\|_{L^\infty}\|b_{x,q}\|_{L^2}\|a_{xx,q}\|_{L^2}\\
&\lesssim \sum_{q\geq -1}\lambda_q^{2s+5}\|a_{q}\|_{L^2}^{2}\|b_{q}\|_{L^2}\\
&\lesssim \sum_{q\geq -1}\left(\lambda_q^{s+2}\|a_{q}\|_{L^2}\right) \left(\lambda_q^{s+1}\|a_{q}\|_{L^2}\lambda_q^{s}\|b_{q}\|_{L^2}\lambda_q^{2-s}\right)\\
&\leq \frac{\mu}{64}\sum_{q\geq -1}\lambda_q^{2s+4}\|a_{q}\|_{L^2}^2+\frac{C}{\mu}\sum_{q\geq -1}\lambda_q^{2s+2}\|a_{q}\|_{L^2}^2\lambda_q^{2s}\|b_{q}\|_{L^2}^2\lambda_q^{2(2-s)}\\
&\leq \frac{\mu}{64}\sum_{q\geq -1}\lambda_q^{2s+4}\|a_{q}\|_{L^2}^2+\frac{C}{\mu}\left(\sum_{q\geq -1}\lambda_q^{2s+2}\|a_{q}\|_{L^2}^2\right)\left(\sum_{q\geq -1}\lambda_q^{2s}\|b_{q}\|_{L^2}^2\right).
\end{split}
\end{equation}

We remark that $K_{111}$  and $K_{113}$ share the same estimates as for $I_{111}$ and $I_{113}$ respectively.

The estimate of $I_{12}$ is given by
\begin{equation}\notag
\begin{split}
|I_{12}|&\leq \sum_{q\geq -1}\sum_{p\sim q}\lambda_q^{2s}\|a_{y,p}\|_{L^2}\|a_{xx,q}\|_{L^2}\|b_{x,\leq p-2}\|_{L^\infty}\\
&\lesssim  \sum_{q\geq -1}\lambda_q^{2s+3}\|a_{q}\|_{L^2}^2\sum_{p\leq q}\lambda_{p}^2\|b_{p}\|_{L^2}\\
&\lesssim  \sum_{q\geq -1}\left(\lambda_q^{s+2}\|a_{q}\|_{L^2}\right) \left(\lambda_q^{s+1}\|a_{q}\|_{L^2}\right)\sum_{p\leq q}\lambda_{p}^s\|b_{p}\|_{L^2}\lambda_p^{2-s}\\
&\leq \frac{\mu}{64}\sum_{q\geq-1} \lambda_q^{2s+4}\|a_{q}\|_{L^2}^2+\frac{C}{\mu} \sum_{q\geq-1} \lambda_q^{2s+2}\|a_q\|_{L^2}^2\left(\sum_{p\leq q} \lambda_p^{s}\|b_p\|_{L^2}\lambda_p^{2-s}\right)^2\\
&\leq \frac{\mu}{64}\sum_{q\geq-1} \lambda_q^{2s+4}\|a_{q}\|_{L^2}^2+\frac{C}{\mu} \left(\sum_{q\geq-1} \lambda_q^{2s+2}\|a_q\|_{L^2}^2\right)\left(\sum_{p\geq-1} \lambda_p^{2s}\|b_p\|_{L^2}^2\right).
\end{split}
\end{equation}
We proceed to estimate of $I_{13}$ as
\begin{equation}\notag
\begin{split}
|I_{13}|&\leq \sum_{p\geq -1}\sum_{q\leq p+2}\lambda_q^{2s}\|\tilde a_{y,p}\|_{L^\infty}\|b_{x,p}\|_{L^2}\|a_{xx,q}\|_{L^2}\\
&\lesssim \sum_{p\geq -1}\lambda_p^3\|a_{p}\|_{L^2}\|b_{p}\|_{L^2}\sum_{q\leq p+2}\lambda_q^{2s+2}\|a_{q}\|_{L^2}\\
&\lesssim \sum_{p\geq -1}\left(\lambda_p^{s+2}\|a_{p}\|_{L^2}\right)\left(\lambda_p^{s}\|b_{p}\|_{L^2}\right)\sum_{q\leq p+2}\lambda_q^{s+1}\|a_{q}\|_{L^2}\lambda_{p-q}^{1-2s}\lambda_q^{2-s}\\
&\leq \frac{\mu}{64}\sum_{p\geq -1}\lambda_p^{2s+4}\|a_{p}\|_{L^2}^2+\frac{C}{\mu}\sum_{p\geq-1}\lambda_p^{2s}\|b_{p}\|_{L^2}^2\left(\sum_{q\leq p+2}\lambda_q^{s+1}\|a_{q}\|_{L^2}\lambda_{p-q}^{1-2s}\right)^2.
\end{split}
\end{equation}
Applying Jensen's inequality yields
\begin{equation}\notag
\begin{split}
&\sum_{p\geq-1}\lambda_p^{2s}\|b_{p}\|_{L^2}^2\left(\sum_{q\leq p+2}\lambda_q^{s+1}\|a_{q}\|_{L^2}\lambda_{p-q}^{1-2s}\right)^2\\
\lesssim& \sum_{p\geq-1}\lambda_p^{2s}\|b_{p}\|_{L^2}^2\left(\sum_{q\leq p+2}\lambda_q^{2s+2}\|a_{q}\|_{L^2}^2\right)\\
\lesssim& \sum_{p\geq-1}\lambda_p^{2s}\|b_{p}\|_{L^2}^2\left(\sum_{q<\infty}\lambda_q^{2s+2}\|a_{q}\|_{L^2}^2\right).
\end{split}
\end{equation}
Hence we have
\begin{equation}\notag
|I_{13}|
\leq \frac{\mu}{64}\sum_{p\geq -1}\lambda_p^{2s+4}\|a_{p}\|_{L^2}^2+\frac{C}{\mu}\sum_{q\geq -1}\lambda_q^{2s+2}\|a_{q}\|_{L^2}^2\sum_{p\geq-1}\lambda_p^{2s}\|b_{p}\|_{L^2}^2.
\end{equation}
The term $K_{12}$ can be estimated similarly as for $I_{111}$
\begin{equation}\notag
\begin{split}
|K_{12}|&\leq \sum_{q\geq -1}\sum_{p\sim q}\lambda_q^{2s}\|a_{y,p}\|_{L^2}\|b_{x,q}\|_{L^2}\|a_{xx,\leq p-2}\|_{L^\infty}\\
&\lesssim  \sum_{q\geq -1}\lambda_q^{2s+2}\|a_{q}\|_{L^2}\|b_{q}\|_{L^2}\sum_{p\leq q}\lambda_{p}^3\|a_{p}\|_{L^2}\\
&\leq \frac{\mu}{64}\sum_{q\geq-1} \lambda_q^{2s+4}\|a_{q}\|_{L^2}^2+\frac{C}{\mu} \sum_{q\geq-1} \lambda_q^{2s+2}\|a_q\|_{L^2}^2\sum_{q\geq-1} \lambda_q^{2s}\|b_q\|_{L^2}^2.
\end{split}
\end{equation}
In the end, we estimate $K_{13}$ as
\begin{equation}\notag
\begin{split}
|K_{13}|&\leq  \sum_{q\geq -1}\sum_{q\leq p+2}\lambda_q^{2s}\|\tilde a_{y,p}\|_{L^2}\| a_{xx,p}\|_{L^2}\|b_{x,q}\|_{L^\infty}\\
&\lesssim  \sum_{q\geq -1}\lambda_q^{2s+2}\|b_{q}\|_{L^2}\sum_{q\leq p+2}\lambda_p^3\|a_{p}\|_{L^2}^2\\
&\lesssim  \sum_{p\geq -1}\lambda_p^3\|a_{p}\|_{L^2}^2\sum_{q\leq p+2}\lambda_q^{2s+2}\|b_{q}\|_{L^2}\\
&\lesssim  \sum_{p\geq -1}\left(\lambda_p^{s+2}\|a_{p}\|_{L^2}\right)\left(\lambda_p^{s+1}\|a_{p}\|_{L^2}\right)\sum_{q\leq p+2}\lambda_q^{s}\|b_{q}\|_{L^2}\lambda_{p-q}^{-2s}\lambda_q^{2-s}\\
&\leq\ \frac{\mu}{64}  \sum_{p\geq -1}\lambda_p^{2s+4}\|a_{p}\|_{L^2}^2+\frac{C}{\mu}\sum_{p\geq -1}\lambda_p^{2s+2}\|a_{p}\|_{L^2}^2\left(\sum_{q\leq p+2}\lambda_q^{s}\|b_{q}\|_{L^2}\lambda_{p-q}^{-2s}\right)^2\\
&\leq\ \frac{\mu}{64}  \sum_{p\geq -1}\lambda_p^{2s+4}\|a_{p}\|_{L^2}^2+\frac{C}{\mu}\sum_{p\geq -1}\lambda_p^{2s+2}\|a_{p}\|_{L^2}^2\left(\sum_{q\leq p+2}\lambda_q^{2s}\|b_{q}\|_{L^2}^2\right)\\
&\leq\ \frac{\mu}{64}  \sum_{p\geq -1}\lambda_p^{2s+4}\|a_{p}\|_{L^2}^2+\frac{C}{\mu}\left(\sum_{p\geq -1}\lambda_p^{2s+2}\|a_{p}\|_{L^2}^2\right)\left(\sum_{q\geq -1}\lambda_q^{2s}\|b_{q}\|_{L^2}^2\right)
\end{split}
\end{equation}
where we have used H\"older's, Bernstein's, Young's and Jensen's inequalities. 

Summarizing the analysis above yields, for $s>2$
\begin{equation}\label{est-I1}
\begin{split}
&|I_1+K_1|\\
\leq&\  \frac{\mu}{8}\sum_{q\geq -1}\lambda_q^{2s+4}\|a_{q}\|_{L^2}^2+\frac{C}{\mu}\left(\sum_{q\geq -1}\lambda_q^{2s+2}\|a_{q}\|_{L^2}^2\right)\left(\sum_{q\geq -1}\lambda_q^{2s}\|b_{q}\|_{L^2}^2\right)\\
\leq&\  \frac{\mu}{8}\sum_{q\geq -1}\lambda_q^{2s+4}\|a_{q}\|_{L^2}^2+\frac{C}{\mu}\left(\sum_{q\geq -1}\lambda_q^{2s+2}\|a_{q}\|_{L^2}^2\right)^2+\frac{C}{\mu}\left(\sum_{q\geq -1}\lambda_q^{2s}\|b_{q}\|_{L^2}^2\right)^2.
\end{split}
\end{equation}

\subsection{Estimate of $I_2+K_2$}
We only provide the decompositions and the use of commutators to show cancellations within $I_2+K_2$. The estimates of the remaining terms can be obtained analogously as for $I_1+K_1$. 
Applying Bony's paraproduct again gives
\begin{equation}\notag
\begin{split}
I_2
=&\sum_{q\geq -1}\sum_{|p-q|\leq 2}\lambda_q^{2s}\int_{\mathbb R^2}(a_{y,\leq p-2}b_{x,p})_q a_{yy,q} \, d\vec x\\
&+\sum_{q\geq -1}\sum_{|p-q|\leq 2}\lambda_q^{2s}\int_{\mathbb R^2}(b_{x,\leq p-2}a_{y,p})_q a_{yy,q} \, d\vec x\\
&+\sum_{q\geq -1}\sum_{p\geq q- 2}\lambda_q^{2s}\int_{\mathbb R^2}(\widetilde a_{y,p}b_{x,p})_q a_{yy,q} \, d\vec x\\
=&: I_{21}+I_{22}+I_{23},
\end{split}
\end{equation}
\begin{equation}\notag
\begin{split}
K_2
=&-\sum_{q\geq-1}\sum_{|p-q|\leq 2} \lambda_q^{2s}\int_{\mathbb R^2} (a_{y,\leq p-2}a_{yy, p})_qb_{x,q}\, d\vec x\\
&-\sum_{q\geq-1}\sum_{|p-q|\leq 2} \lambda_q^{2s}\int_{\mathbb R^2} (a_{yy,\leq p-2}a_{y, p})_qb_{x,q}\, d\vec x\\
&-\sum_{q\geq-1}\sum_{p\geq q- 2} \lambda_q^{2s}\int_{\mathbb R^2} (\widetilde a_{y,p}a_{yy, p})_qb_{x,q}\, d\vec x\\
=&: K_{21}+K_{22}+K_{23}.
\end{split}
\end{equation}
Applying commutators we further have
\begin{equation}\notag
\begin{split}
I_{21}=&\sum_{q\geq -1}\sum_{|p-q|\leq 2}\lambda_q^{2s} \int_{\mathbb R^2} [\Delta_q, a_{y,\leq p-2}\partial_x] b_p a_{yy,q} \, d\vec x\\
&+\sum_{q\geq -1}\sum_{|p-q|\leq 2}\lambda_q^{2s} \int_{\mathbb R^2} a_{y,\leq q-2}(b_{x,p})_q a_{yy,q} \, d\vec x\\
&+\sum_{q\geq -1}\sum_{|p-q|\leq 2}\lambda_q^{2s} \int_{\mathbb R^2} \left(a_{y,\leq q-2}-a_{y,\leq p-2}\right)(b_{x,p})_q a_{yy,q} \, d\vec x\\
=&:I_{211}+I_{212}+I_{213},
\end{split}
\end{equation}
\begin{equation}\notag
\begin{split}
K_{21}=&-\sum_{q\geq -1}\sum_{|p-q|\leq 2}\lambda_q^{2s} \int_{\mathbb R^2} [\Delta_q, a_{y,\leq p-2}\partial_y] a_{y,p} b_{x,q} \, d\vec x\\
&-\sum_{q\geq -1}\sum_{|p-q|\leq 2}\lambda_q^{2s} \int_{\mathbb R^2} a_{y,\leq q-2}(a_{yy,p})_q b_{x,q} \, d\vec x\\
&-\sum_{q\geq -1}\sum_{|p-q|\leq 2}\lambda_q^{2s} \int_{\mathbb R^2} \left(a_{y,\leq p-2}-a_{y,\leq q-2}\right)(a_{yy,p})_q b_{x,q} \, d\vec x\\
=&:K_{211}+K_{212}+K_{213}.
\end{split}
\end{equation}
As before, we observe the cancellation 
\[I_{212}+K_{212}=0\]
thanks to $\sum_{|p-q|\leq 2}(b_{x,p})_q=b_{x,q}$ and $\sum_{|p-q|\leq 2}(a_{yy,p})_q=a_{yy,q}$ for any fixed $q\geq -1$.

Noticing that the only difference between $I_2+K_2$ and $I_1+K_1$ is that, the former contains $a_{yy}$ while the later contains $a_{xx}$. Such difference does not affect the estimates. Hence we claim 
\begin{equation}\label{est-I2}
\begin{split}
&|I_2+K_2|\\
\leq&\  \frac{\mu}{8}\sum_{q\geq -1}\lambda_q^{2s+4}\|a_{q}\|_{L^2}^2+\frac{C}{\mu}\left(\sum_{q\geq -1}\lambda_q^{2s+2}\|a_{q}\|_{L^2}^2\right)^2+\frac{C}{\mu}\left(\sum_{q\geq -1}\lambda_q^{2s}\|b_{q}\|_{L^2}^2\right)^2.
\end{split}
\end{equation}

\subsection{The a priori estimate}
Applying the fact that $I_3+K_3$ and $I_4+K_4$ share the same estimates of $I_2+K_2$ and $I_1+K_1$ respectively, combining \eqref{priori-3}, \eqref{est-I1} and \eqref{est-I2} gives
\begin{equation}\label{priori-4}
\begin{split}
&\frac{d}{dt} \sum_{q\geq-1}\left( \lambda_q^{2s+2}\| a_q\|_{L^2}^2+ \lambda_q^{2s}\|b_q\|_{L^2}^2\right)+\mu \sum_{q\geq-1} \lambda_q^{2s+4}\|a_q\|_{L^2}^2\\
\leq &\ \frac{C}{\mu}\left(\sum_{q\geq -1}\lambda_q^{2s+2}\|a_{q}\|_{L^2}^2\right)^2+\frac{C}{\mu}\left(\sum_{q\geq -1}\lambda_q^{2s}\|b_{q}\|_{L^2}^2\right)^2.
\end{split}
\end{equation}
It follows from Gr\"onwall's inequality that there exists a time $T>0$ such that on $[0,T)$
\begin{equation}\notag
\sum_{q\geq-1}\left( \lambda_q^{2s+2}\| a_q(t)\|_{L^2}^2+ \lambda_q^{2s}\|b_q(t)\|_{L^2}^2\right)\lesssim \|a_0\|_{H^{s+1}}^2+\|b_0\|_{H^{s}}^2.
\end{equation}
Combining the basic energy estimate in Subsection \ref{sec-basic}, it completes the proof the a priori estimate for $(a(t), b(t))\in H^{s+1}\times H^s$.

\bigskip

\section{Proof of Theorem \ref{thm} for $\mu=0$ and $\nu>0$}
\label{sec-proof2}

As argued in Section \ref{sec-proof1}, it is sufficient to show the a priori estimate for $(a(t), b(t))\in H^{s+1}\times H^s$ with $s>2$.
We also note it is enough to provide details for the estimate of $I_1+K_1$. Assuming $\mu=0$ and $\nu>0$, we estimate the non-zero terms of $I_1+K_1$.

As before, applying H\"older's, Bernstein's and Young's inequalities yields
\begin{equation}\notag
\begin{split}
|I_{111}|&\leq \sum_{q\geq-1}\sum_{|p-q|\leq 2} \lambda_q^{2s} \|a_{xy, \leq p-2}\|_{L^\infty}\|b_p\|_{L^2}\|a_{xx,q}\|_{L^2}\\
&\lesssim \sum_{q\geq-1} \lambda_q^{2s+2}\|a_{q}\|_{L^2}\|b_q\|_{L^2}\sum_{p\leq q}\lambda_p^{3} \|a_{p}\|_{L^2}\\
&\lesssim \sum_{q\geq-1} \left(\lambda_q^{s+1}\|b_{q}\|_{L^2}\right)\left(\lambda_q^{s+1}\|a_q\|_{L^2}\right)\sum_{p\leq q}\lambda_p^{s+1} \|a_{p}\|_{L^2}\lambda_{p}^{2-s}\\
&\leq \frac{\nu}{64}\sum_{q\geq-1} \lambda_q^{2s+2}\|b_{q}\|_{L^2}^2+\frac{C}{\nu} \sum_{q\geq-1} \lambda_q^{2s+2}\|a_q\|_{L^2}^2\left(\sum_{p\leq q}\lambda_p^{s+1} \|a_{p}\|_{L^2}\lambda_{p}^{2-s}\right)^2.
\end{split}
\end{equation}
Applying H\"older's inequality for sequence gives, for $s>2$
\begin{equation}\notag
\begin{split}
|I_{111}|
&\leq \frac{\nu}{64}\sum_{q\geq-1} \lambda_q^{2s+2}\|b_{q}\|_{L^2}^2+\frac{C}{\nu}\left(\sum_{p<\infty}\lambda_p^{s+1} \|a_{p}\|_{L^2}\lambda_{p}^{2-s}\right)^2 \sum_{q\geq-1} \lambda_q^{2s+2}\|a_q\|_{L^2}^2\\
&\leq \frac{\nu}{64}\sum_{q\geq-1} \lambda_q^{2s+2}\|b_{q}\|_{L^2}^2+\frac{C}{\nu}\sum_{p<\infty}\lambda_{p}^{2(2-s)}\sum_{p<\infty}\lambda_p^{2s+2} \|a_{p}\|_{L^2}^2\sum_{q\geq-1} \lambda_q^{2s+2}\|a_q\|_{L^2}^2\\
&\leq \frac{\nu}{64}\sum_{q\geq-1} \lambda_q^{2s+2}\|b_{q}\|_{L^2}^2+\frac{C}{\nu} \left(\sum_{q\geq-1} \lambda_q^{2s+2}\|a_q\|_{L^2}^2\right)^2.
\end{split}
\end{equation}
The estimate of $I_{113}$ follows similarly
\begin{equation}\notag
\begin{split}
|I_{113}|&\lesssim \sum_{q\geq -1}\lambda_q^{2s}\|a_{y,q}\|_{L^\infty}\|b_{x,q}\|_{L^2}\|a_{xx,q}\|_{L^2}\\
&\lesssim \sum_{q\geq -1}\lambda_q^{2s+5}\|a_{q}\|_{L^2}^{2}\|b_{q}\|_{L^2}\\
&\lesssim \sum_{q\geq -1}\left(\lambda_q^{s+1}\|b_{q}\|_{L^2}\right) \left(\lambda_q^{2s+2}\|a_{q}\|_{L^2}^{2}\lambda_q^{2-s}\right)\\
&\leq \frac{\nu}{64}\sum_{q\geq -1}\lambda_q^{2s+2}\|b_{q}\|_{L^2}^2+\frac{C}{\nu}\sum_{q\geq -1}\lambda_q^{4s+4}\|a_{q}\|_{L^2}^{4}\lambda_q^{2(2-s)}.
\end{split}
\end{equation}
Using the fact $s>2$ again, it follows
\begin{equation}\notag
\begin{split}
\sum_{q\geq -1}\lambda_q^{4s+4}\|a_{q}\|_{L^2}^{4}\lambda_q^{2(2-s)}
&\leq  \left( \sum_{q\geq -1}\lambda_q^{s+1}\|a_{q}\|_{L^2}\lambda_q^{\frac12(2-s)}\right)^4\\
&\leq  \left( \sum_{q\geq -1}\lambda_q^{2s+2}\|a_{q}\|_{L^2}^2\right)^2.
\end{split}
\end{equation}
Therefore we obtain
\begin{equation}\notag
\begin{split}
|I_{113}|
&\leq \frac{\nu}{64}\sum_{q\geq -1}\lambda_q^{2s+2}\|b_{q}\|_{L^2}^2+\frac{C}{\nu} \left( \sum_{q\geq -1}\lambda_q^{2s+2}\|a_{q}\|_{L^2}^2\right)^2.
\end{split}
\end{equation}

As before we note $K_{111}$  and $K_{113}$ can be respectively estimated analogously as for $I_{111}$ and $I_{113}$.

The term $I_{12}$ is the most difficult term to be estimated since more derivatives are on the high modes. To over come the obstacle, we need to apply commutator to move one derivative to the low modes part. We thus rearrange $I_{12}$ as 
\begin{equation}\notag
\begin{split}
I_{12}=&\sum_{q\geq -1}\sum_{|p-q|\leq 2}\lambda_q^{2s}\int_{\mathbb T^2}[\Delta_q, b_{x,\leq p-2}]\partial_ya_{p} a_{xx,q}\, dxdy\\
&+\sum_{q\geq -1}\sum_{|p-q|\leq 2}\lambda_q^{2s}\int_{\mathbb T^2}b_{x,\leq q-2}\Delta_q(a_{y,p}) a_{xx,q}\, dxdy\\
&+\sum_{q\geq -1}\sum_{|p-q|\leq 2}\lambda_q^{2s}\int_{\mathbb T^2}((b_{x,\leq p-2}-b_{x,\leq q-2})\Delta_q(a_{y,p}) a_{xx,q}\, dxdy\\
=&: I_{121}+ I_{122}+ I_{123}.
\end{split}
\end{equation}
Applying H\"older's inequality, the commutator estimate in Lemma \ref{le-comm} and Young's inequality gives
\begin{equation}\notag
\begin{split}
I_{121}\leq &\sum_{q\geq -1}\sum_{|p-q|\leq 2}\lambda_q^{2s}\|[\Delta_q, b_{x,\leq p-2}]\partial_ya_{p} \|_{L^2}\|a_{xx,q}\|_{L^2}\\
\lesssim &\sum_{q\geq -1}\lambda_q^{2s+2} \|a_{q}\|_{L^2} \sum_{|p-q|\leq 2}\|b_{xy,\leq p-2}\|_{L^\infty}\|a_{p}\|_{L^2}\\
\lesssim &\sum_{q\geq -1}\lambda_q^{2s+2} \|a_{q}\|_{L^2}^2 \sum_{p\leq q+2}\lambda_p^3\|b_{p}\|_{L^2}\\
\lesssim &\sum_{p\geq -1}\left(\lambda_p^{s+1}\|b_{p}\|_{L^2}\right)\lambda_p^{2-s}\sum_{q\geq p-2}\lambda_q^{2s+2} \|a_{q}\|_{L^2}^2\\
\leq &\ \frac{\nu}{64}\sum_{p\geq -1}\lambda_p^{2s+2}\|b_{p}\|_{L^2}^2+\frac{C}{\nu}\sum_{p\geq -1}\lambda_p^{2(2-s)}\left(\sum_{q\geq p-2}\lambda_q^{2s+2} \|a_{q}\|_{L^2}^2\right)^2\\
\leq &\ \frac{\nu}{64}\sum_{p\geq -1}\lambda_p^{2s+2}\|b_{p}\|_{L^2}^2+\frac{C}{\nu}\left(\sum_{q\geq -1}\lambda_q^{2s+2} \|a_{q}\|_{L^2}^2\right)^2.
\end{split}
\end{equation}
Applying the fact $\sum_{|p-q|\leq 2}\Delta_q(a_{y,p})=\Delta_q a_y$ for any fixed $q\geq -1$ and integration by parts to $I_{122}$ we deduce
\begin{equation}\notag
\begin{split}
I_{122}
=&\sum_{q\geq -1}\sum_{|p-q|\leq 2}\lambda_q^{2s}\int_{\mathbb T^2}b_{x,\leq q-2}a_{y,q} a_{xx,q}\, dxdy\\
=&-\sum_{q\geq -1}\sum_{|p-q|\leq 2}\lambda_q^{2s}\int_{\mathbb T^2}b_{xx,\leq q-2}a_{y,q} a_{x,q}\, dxdy\\
&-\sum_{q\geq -1}\sum_{|p-q|\leq 2}\lambda_q^{2s}\int_{\mathbb T^2}b_{x,\leq q-2}a_{xy,q} a_{x,q}\, dxdy.
\end{split}
\end{equation}
Noticing $a_{xy,q} a_{x,q}=\frac12\partial_y (a_{x,q})^2$, we apply integration by parts again
\begin{equation}\notag
\begin{split}
I_{122}
=&-\sum_{q\geq -1}\sum_{|p-q|\leq 2}\lambda_q^{2s}\int_{\mathbb T^2}b_{xx,\leq q-2}a_{y,q} a_{x,q}\, dxdy\\
&-\frac12\sum_{q\geq -1}\sum_{|p-q|\leq 2}\lambda_q^{2s}\int_{\mathbb T^2}b_{x,\leq q-2}\partial_y (a_{x,q})^2\, dxdy\\
=&-\sum_{q\geq -1}\sum_{|p-q|\leq 2}\lambda_q^{2s}\int_{\mathbb T^2}b_{xx,\leq q-2}a_{y,q} a_{x,q}\, dxdy\\
&+\frac12\sum_{q\geq -1}\sum_{|p-q|\leq 2}\lambda_q^{2s}\int_{\mathbb T^2}b_{xy,\leq q-2}(a_{x,q})^2\, dxdy.
\end{split}
\end{equation}
Then as usual, we apply H\"older's inequality, Bernstein's inequality
\begin{equation}\notag
\begin{split}
|I_{122}|
\lesssim &\sum_{q\geq -1}\sum_{|p-q|\leq 2}\lambda_q^{2s}\|b_{xx,\leq q-2}\|_{L^\infty}\|a_{y,q}\|_{L^2}\|a_{x,q}\|_{L^2}\\
&+\sum_{q\geq -1}\sum_{|p-q|\leq 2}\lambda_q^{2s}\|b_{xy,\leq q-2}\|_{L^\infty}\|a_{x,q}\|_{L^2}^2\\
\lesssim &\sum_{q\geq -1}\lambda_q^{2s+2} \|a_{q}\|_{L^2}^2 \sum_{p\leq q+2}\lambda_p^3\|b_{p}\|_{L^2}\\
\leq &\ \frac{\nu}{64}\sum_{p\geq -1}\lambda_p^{2s+2}\|b_{p}\|_{L^2}^2+\frac{C}{\nu}\left(\sum_{q\geq -1}\lambda_q^{2s+2} \|a_{q}\|_{L^2}^2\right)^2
\end{split}
\end{equation}
where the last step is obtained similarly as for $I_{121}$. 

Realizing that there is only a finite number of modes in the term $b_{x,\leq p-2}-b_{x,\leq q-2}$ for $|p-q|\leq 2$, we estimate $I_{123}$ as
\begin{equation}\notag
\begin{split}
|I_{123}|
\leq &\sum_{q\geq -1}\sum_{|p-q|\leq 2}\lambda_q^{2s}\|b_{x,\leq p-2}-b_{x,\leq q-2}\|_{L^\infty}\|a_{y,p}\|_{L^2}\|a_{xx,q}\|_{L^2}\\
\lesssim &\sum_{q\geq -1}\lambda_q^{2s+3} \|a_{q}\|_{L^2}^2 \sum_{|p-q|\leq 2}\lambda_p\|b_{p}\|_{L^\infty}\\
\lesssim &\sum_{q\geq -1}\lambda_q^{2s+5} \|a_{q}\|_{L^2}^2 \|b_{q}\|_{L^2}\\
\leq &\ \frac{\nu}{64}\sum_{p\geq -1}\lambda_p^{2s+2}\|b_{p}\|_{L^2}^2+\frac{C}{\nu}\left(\sum_{q\geq -1}\lambda_q^{2s+2} \|a_{q}\|_{L^2}^2\right)^2
\end{split}
\end{equation}
where we obtain the last step analogously as for $I_{113}$.

The estimate of $I_{13}$ is less complicated, given by
\begin{equation}\notag
\begin{split}
|I_{13}|&\leq \sum_{p\geq -1}\sum_{q\leq p+2}\lambda_q^{2s}\|\tilde a_{y,p}\|_{L^\infty}\|b_{x,p}\|_{L^2}\|a_{xx,q}\|_{L^2}\\
&\lesssim \sum_{p\geq -1}\lambda_p^3\|a_{p}\|_{L^2}\|b_{p}\|_{L^2}\sum_{q\leq p+2}\lambda_q^{2s+2}\|a_{q}\|_{L^2}\\
&\lesssim \sum_{p\geq -1}\left(\lambda_p^{s+1}\|b_{p}\|_{L^2}\right)\left(\lambda_p^{s+1}\|a_{p}\|_{L^2}\right)\sum_{q\leq p+2}\lambda_q^{s+1}\|a_{q}\|_{L^2}\lambda_{p-q}^{1-2s}\lambda_q^{2-s}\\
&\leq \frac{\nu}{64}\sum_{p\geq -1}\lambda_p^{2s+2}\|b_{p}\|_{L^2}^2+\frac{C}{\nu}\sum_{q\geq-1}\lambda_p^{2s+2}\|a_{p}\|_{L^2}^2\left(\sum_{q\leq p+2}\lambda_q^{s+1}\|a_{q}\|_{L^2}\lambda_{p-q}^{1-2s}\right)^2.
\end{split}
\end{equation}
We continue with Jensen's inequality 
\begin{equation}\notag
\begin{split}
&\sum_{q\geq-1}\lambda_p^{2s+2}\|a_{p}\|_{L^2}^2\left(\sum_{q\leq p+2}\lambda_q^{s+1}\|a_{q}\|_{L^2}\lambda_{p-q}^{1-2s}\right)^2\\
\lesssim& \sum_{q\geq-1}\lambda_p^{2s+2}\|a_{p}\|_{L^2}^2\left(\sum_{q\leq p+2}\lambda_q^{2s+2}\|a_{q}\|_{L^2}^2\right)\\
\lesssim& \sum_{q\geq-1}\lambda_p^{2s+2}\|a_{p}\|_{L^2}^2\left(\sum_{q<\infty}\lambda_q^{2s+2}\|a_{q}\|_{L^2}^2\right).
\end{split}
\end{equation}
Combining the last two estimates gives
\begin{equation}\notag
|I_{13}|
\leq \frac{\nu}{64}\sum_{p\geq -1}\lambda_p^{2s+2}\|b_{p}\|_{L^2}^2+\frac{C}{\nu} \left(\sum_{q\geq -1}\lambda_q^{2s+2}\|a_{q}\|_{L^2}^2\right)^2.
\end{equation}

The estimate of $K_{12}$ is also simple since it can be treated similarly as for $I_{111}$
\begin{equation}\notag
\begin{split}
|K_{12}|&\leq \sum_{q\geq -1}\sum_{p\sim q}\lambda_q^{2s}\|a_{y,p}\|_{L^2}\|b_{x,q}\|_{L^2}\|a_{xx,\leq p-2}\|_{L^\infty}\\
&\lesssim  \sum_{q\geq -1}\lambda_q^{2s+2}\|a_{q}\|_{L^2}\|b_{q}\|_{L^2}\sum_{p\leq q}\lambda_{p}^3\|a_{p}\|_{L^2}\\
&\leq \frac{\nu}{64}\sum_{q\geq-1} \lambda_q^{2s+2}\|b_{q}\|_{L^2}^2+\frac{C}{\nu} \left(\sum_{q\geq-1} \lambda_q^{2s+2}\|a_q\|_{L^2}^2\right)^2.
\end{split}
\end{equation}

Applying H\"older's, Bernstein's and Young's inequalities to the last term in $I_1+K_1$ gives
\begin{equation}\notag
\begin{split}
|K_{13}|&\leq  \sum_{q\geq -1}\sum_{q\leq p+2}\lambda_q^{2s}\|\tilde a_{y,p}\|_{L^2}\| a_{xx,p}\|_{L^2}\|b_{x,q}\|_{L^\infty}\\
&\lesssim  \sum_{q\geq -1}\lambda_q^{2s+2}\|b_{q}\|_{L^2}\sum_{q\leq p+2}\lambda_p^3\|a_{p}\|_{L^2}^2\\
&\lesssim  \sum_{q\geq -1}\left(\lambda_q^{s+1}\|b_{q}\|_{L^2}\right) \lambda_q^{2-s}\sum_{q\leq p+2}\lambda_p^{2s+2}\|a_{p}\|_{L^2}^2\lambda_{p-q}^{1-2s}\\
&\leq\ \frac{\nu}{64}  \sum_{q\geq -1}\lambda_q^{2s+2}\|b_{q}\|_{L^2}^2+\frac{C}{\nu}\sum_{q\geq -1} \lambda_q^{2(2-s)}\left(\sum_{q\leq p+2}\lambda_p^{2s+2}\|a_{p}\|_{L^2}^2\lambda_{p-q}^{1-2s}\right)^2
\end{split}
\end{equation}
where we used $s>2$. 
It the follows from Jensen's inequality
\begin{equation}\notag
\begin{split}
|K_{13}|
&\leq\ \frac{\nu}{64}  \sum_{q\geq -1}\lambda_q^{2s+2}\|b_{q}\|_{L^2}^2+\frac{C}{\nu}\sum_{q\geq -1} \lambda_q^{2(2-s)}\sum_{q\leq p+2}\lambda_p^{2(2s+2)}\|a_{p}\|_{L^2}^4\\
&\leq\ \frac{\nu}{64}  \sum_{q\geq -1}\lambda_q^{2s+2}\|b_{q}\|_{L^2}^2+\frac{C}{\nu}\sum_{p<\infty} \lambda_p^{2(2s+2)}\|a_{p}\|_{L^2}^4\\
&\leq\ \frac{\nu}{64}  \sum_{q\geq -1}\lambda_q^{2s+2}\|b_{q}\|_{L^2}^2+\frac{C}{\nu}\left(\sum_{p\geq -1} \lambda_p^{2s+2}\|a_{p}\|_{L^2}^2\right)^2.
\end{split}
\end{equation}

Therefore the analysis above implies 
\begin{equation}\label{est-I1-nu}
|I_1+K_1|\leq \frac{\nu}{8}  \sum_{q\geq -1}\lambda_q^{2s+2}\|b_{q}\|_{L^2}^2+\frac{C}{\nu}\left(\sum_{p\geq -1} \lambda_p^{2s+2}\|a_{p}\|_{L^2}^2\right)^2,
\end{equation}
and $I_2+K_2$, $I_3+K_3$, $I_4+K_4$ also satisfy this estimate. Combining \eqref{priori-3} with $\mu=0$ and $\nu>0$ and \eqref{est-I1-nu} yields
\begin{equation}\label{priori-5}
\begin{split}
&\frac{d}{dt} \sum_{q\geq-1}\left( \lambda_q^{2s+2}\| a_q\|_{L^2}^2+ \lambda_q^{2s}\|b_q\|_{L^2}^2\right)+\nu \sum_{q\geq-1} \lambda_q^{2s+2}\|b_q\|_{L^2}^2\\
\leq &\ \frac{C}{\nu}\left(\sum_{q\geq -1}\lambda_q^{2s+2}\|a_{q}\|_{L^2}^2\right)^2.
\end{split}
\end{equation}
Applying Gr\"onwall's inequality again to \eqref{priori-5} indicates that 
\begin{equation}\notag
\sum_{q\geq-1}\left( \lambda_q^{2s+2}\| a_q(t)\|_{L^2}^2+ \lambda_q^{2s}\|b_q(t)\|_{L^2}^2\right)\lesssim \|a_0\|_{H^{s+1}}^2+\|b_0\|_{H^{s}}^2, \ \ t\in[0,T)
\end{equation}
for some $T>0$.




\bigskip



\begin{thebibliography}{XX}


































\bibitem{Dai-emhd-25}
M. Dai.
\newblock {\em Blowup for the forced electron MHD}.
\newblock  arXiv: 2503.14777, 2025.

\bibitem{Dai-emhd-23}
M. Dai.
\newblock {\em Global existence of 2D electron MHD near a steady state}.
\newblock  arXiv: 2306.13036, 2023.

\bibitem{Dai-emhd-24}
M. Dai.
\newblock {\em Ill-posedness of $2\frac12$D electron MHD}.
\newblock  arXiv: 2411.00120, 2024.




\bibitem{Dai-emhd-18}
M. Dai.
\newblock {\em Non-unique weak solutions in Leray-Hopf class of the 3D Hall-MHD system}.
\newblock SIAM Journal of Mathematical Analysis, 53(5): 5979--6016, 2021.

\bibitem{Dai-emhd-1d}
M. Dai.
\newblock {\em Well-posedness and blowup of 1D electron magnetohydrodynamics}.
\newblock arXiv:2503.06383, 2025.





\bibitem{Dai-Oh}
M. Dai and S.J. Oh.
\newblock {\em Beale--Kato--Majda-type continuation criteria for Hall- and electron-magnetohydrodynamics}.
\newblock arXiv:2407.04314, 2024.

\bibitem{Dai-W}
M. Dai and C. Wu.
\newblock {\em Dissipation wavenumber and regularity for electron magnetohydrodynamics}.
\newblock Journal of Differential Equations, Vol. 376: 655--681, 2023.



































\bibitem{JO1}
I. Jeong and S. Oh.
\newblock {\em On illposedness of the Hall and electron magnetohydrodynamic equations without resistivity on the whole space}.
\newblock arXiv: 2404.13790, 2024.

\bibitem{JO2}
I. Jeong and S. Oh.
\newblock {\em On the Cauchy problem for the Hall and electron magnetohydrodynamic equations without resistivity I: illposedness near degenerate stationary solutions}.
\newblock Annals of PDE, vol.8, no.15, 2022.

\bibitem{JO3}
I. Jeong and S. Oh.
\newblock {\em Wellposedness of the electron MHD without resistivity for large perturbations of the uniform magnetic field}.
\newblock arXiv: 2402.06278, 2024.








































\end{thebibliography}
\end{document}